\documentclass[11pt]{article}

\usepackage{amsfonts,epsf,amsmath,amssymb,here}
\usepackage[numbers]{natbib}
\usepackage{epsfig}
\usepackage{latexsym,multirow,rotating}
\usepackage{pgf,tikz,subfigure}
\usetikzlibrary{decorations.markings}

\newtheorem{theorem}{\bf Theorem}[section]

\newtheorem{proposition}[theorem]{\bf Proposition}

\newtheorem{problem}[theorem]{\bf Problem}

\newcommand{\proof}{\noindent{\bf Proof.\ }}
\newcommand{\qed}{\hfill $\square$ \bigskip}

\begin{document}
\baselineskip=0.24in

\vspace*{40mm}

\begin{center}
{\LARGE \bf\sc Strong Traces Model of Self-Assembly Polypeptide Structures}
\bigskip \bigskip

{\large \sc Ga\v sper Fijav\v z}

\smallskip
{\em Faculty of Computer and Information Science, University of Ljubljana, Slovenia} \\
e-mail: {\tt gasper.fijavz@fri.uni-lj.si}

\bigskip

{\large \sc Toma\v z Pisanski}

\smallskip
{\em Faculty of Mathematics and Physics, University of Ljubljana, Slovenia} \\
e-mail: {\tt tomaz.pisanski@fmf.uni-lj.si}

\bigskip

{\large \sc Jernej Rus}\\
\smallskip
{\em Faculty of Mathematics and Physics, University of Ljubljana, Slovenia} \\
e-mail: {\tt jernej.rus@gmail.com}

\bigskip\medskip

%(Received July 19, 2013)

\end{center}

% \footnotetext[3] {Corresponding author. If possible, send your
% correspondence via e-mail.}
\noindent
{\bf Abstract}

\vspace{3mm}\noindent

A novel self-assembly strategy for polypeptide nanostructure design was presented in [Design of a single-chain polypeptide tetrahedron assembled from coiled-coil segments, Nature Chemical Biology 9 (2013) 362--366]. The first mathematical model (polypeptide nanostructure can naturally be presented as a skeleton graph of a polyhedron) from [Stable traces as a model for self-assembly of polypeptide nanoscale polyhedrons, MATCH Commun.\ Math.\ Comput.\ Chem. 70 (2013) 317--330] introduced stable traces as the appropriate mathematical description, yet we find them deficient in modeling graphs with either very small ($\le 2$) or large ($\ge 6$) degree vertices.
We introduce \emph{strong traces} which remedy both of the above mentioned drawbacks.  We show that \emph{every} connected graph admits a strong trace by studying a connection between strong traces and graph embeddings. Further we also characterize graphs which admit \emph{parallel (resp. antiparallel) strong traces}.

\vspace{5mm}

\baselineskip=0.24in

%\noindent {\bf Key words}: nanostructure design; strong trace; $d$-repetition; $d$-stable trace; $1$-face embedding; parallel strong trace; antiparallel strong trace
%
%\medskip\noindent
%{\bf AMS subject classification (2010)}: 05C10, 05C45

%%%%%%%%%%%%%%%%%
\section{Introduction}
\label{sec:introd}

Recently Gradi\v sar et.\ al~\cite{gr-2013} presented a novel self-assembly strategy for polypeptide nanostructure design that represents a significant development in biotechnology. The main success of their research is a construction of a polypeptide self-assembling tetrahedron by concatenating $12$ coiled-coil-forming segments in a prescribed order. More precisely, a single polypeptide chain consisting of 12 segments was routed through $6$ edges of the tetrahedron in such a way that every edge was traversed exactly twice. In this way $6$ coiled-coil dimers were created and interlocked into a stable tetrahedral structure. 

A polyhedron $P$ which is composed from a single polymer chain can be naturally represented with a graph $G(P)$ of the polyhedron. As in the self-assembly process every edge of $G(P)$ corresponds to a coiled-coil dimer, exactly two segments are associated with every edge of $G(P)$. 

The first mathematical model was introduced in~\cite{kl-2013}, where the authors have shown that a polyhedral graph $P$ can be realized by interlocking pairs of polypeptide chains if its corresponding graph $G(P)$ contains a stable trace (to be defined later).

We find that the mathematical model introduced in~\cite{kl-2013} has two important deficiencies:
\begin{enumerate}
\item[(1)] it does not account for vertices of degree $\le 2$, and
\item[(2)] it does not successfully model vertices of degree $\ge 6$.
\end{enumerate}

The model proposed in this paper settles the above issues. On one hand it successfully extends to graphs with smaller vertex degrees. Even if for every polyhedron $P$ its graph $G(P)$ has minimum degree $\ge 3$, the model should also include graphs of smaller degrees. It is plausible that a quest may require constructions of polypeptide nanostructures with reactive parts being pendant to the main body of the polyhedral structure.

Now (2) touches the question of defining vertices in our structure. An edge in a graph is defined via identifying pairs of segments along a walk $W$. A vertex on the other hand is only defined implicitly: pairs of segments/edges that lie consecutively on this walk should meet in a common endvertex. 

In the case where no vertex in $G$ has degree $\ge 6$ the procedure --- (i) find a stable trace $W$ in $G$ and (ii) identify pairs of edges along $W$ and fold the resulting structure into a graph --- shall produce the initial graph $G$. However if $G$ has a vertex of degree larger than $6$ this may not be the case. A stable trace in $G$ may fold to a graph different from $G$, as a vertex of degree $\ge 6$ may indeed split into a collection of independent vertices of degree $\ge 3$, see also Fig.~\ref{fig:repetition}.

A \emph{strong trace} in a graph, our key object (to be defined later), successfully resolves both above issues. In one sweep we can model graphs with vertices of both low and high degrees. What is more, strong traces admit a natural connection to embeddings of graphs in higher surfaces.

Our main results state that every connected graph admits a strong trace, and can therefore be correctly realized by folding its strong trace by edge identifications (Theorems~\ref{thm:strong} and \ref{thm:realize}). 

In what follows we use Section~\ref{sec:double} to describe some basic and necessary tools from graph theory. In Section~\ref{sec:main} we connect strong traces and embeddings of graphs and ultimately prove our main results. In Sections~\ref{sec:anti} and~\ref{sec:parallel} we generalize two additional concepts from~\cite{kl-2013,rus-2013} --- antiparallel strong traces, parallel strong traces and also parallel $d$-stable traces.

%%%%%%%%%%%%%%%%%%
\section{Double traces}
\label{sec:double}

All graphs considered in this paper will be connected and finite. We denote the degree of a vertex $v$ by $d_G(v)$. The minimum and the maximum degree of $G$ will be denoted by $\delta(G)$ and $\Delta(G)$, respectively. 

If $v$ is a vertex then $N(v)$ denotes set of vertices adjacent to $v$, and $E(v)$ is the set of edges incident with $v$. If $A$ is a set of vertices then $E(v,A)$ denotes the collection of edges incident with both $v$ and a vertex from $A$.

A \emph{walk} in $G$ is an alternating sequence 
\begin{equation}
W=v_0 e_1 v_1 \ldots v_{\ell-1} e_\ell v_\ell,
\label{eq:walk}
\end{equation}
so that for every $i=1,\ldots,k$ $e_i$ is an edge between vertices $v_{i-1}$ and $v_i$. We say that $W$ \emph{passes through} or \emph{traverses} edges and vertices contained in the sequence~\eqref{eq:walk}. The length of a walk is the number of edges in the sequence, and we call $v_0$ and $v_\ell$ the \emph{endvertices} of $W$. A walk is \emph{closed} if its endvertices coincide. 

An \emph{Euler tour} in $G$ is a closed walk which traverses every edge of $G$ exactly once. $G$ is an \emph{Eulerian graph} if it admits an Euler tour. The fundamental Euler's theorem asserts that a (connected) graph $G$ is Eulerian if and only if all of its vertices are of even degree. 

A \emph{double trace} in $G$ is a closed walk which traverses every edge of $G$ exactly twice. Next result essentially goes back to Euler and was since observed by various authors. 

\begin{proposition}
\label{prop:double}
Every connected graph $G$ has a double trace.
\end{proposition}

Let $W$ be a double trace of length $\ell$, and let $N \subseteq N(v)$. We say that $W$ has an \emph{$N$-repetition at $v$} if the following implication holds:
\begin{equation}
\text{\emph{for every $i \in \{0,\ldots,\ell-1\}$: if $v=v_i$ then $v_{i+1} \in N$ if and only if $v_{i-1} \in N$.}}
\label{eq:repetition}
\end{equation}
Intuitively $W$ has an $N$-repetition at $v$ if whenever $W$ visits $v$ coming from a vertex in $N$ it also returns to a vertex of $N$. Let us also note that we treat a double trace as a closed walk taking indices in~\eqref{eq:repetition} modulo $\ell$. This implies that $v_1$ is the vertex immediately following $v_\ell$.

An $N$-repetition (at $v$) is a \emph{$d$-repetition} if $|N|=d$, and a $d$-repetition will also be called a repetition \emph{of order $d$}. An $N$-repetition at $v$ is \emph{trivial} if $N=\emptyset$ or $N=N(v)$. Clearly if $W$ has an $N$-repetition at $v$, then it also has an $N(v)\setminus N$-repetition at $v$. We shall call this observation \emph{symmetry of repetitions}.
 
In~\cite{kl-2013} a $1$-repetition at $v$ was named a \emph{retracing (at a vertex $v$)}, 
and a $2$-repetition at $v$ was denoted as a \emph{repetition at a vertex $v$}. Note that in this paper a \emph{repetition at $v$} can be of order different than $2$.

We call a double trace without nontrivial repetitions of order $<d$ a \emph{$d$-stable trace}, extending the terms used in~\cite{kl-2013} (where a $1$-stable trace was named a \emph{proper trace} and the term \emph{stable trace} was used to name $2$-stable traces). Graphs which admit $1$-stable traces were independently characterized by Sabidussi~\cite{sa-1977} and later by Eggleton and Skilton~\cite{eg-1984}. Graphs admitting $2$-stable traces were recently characterized in~\cite{kl-2013}: 

\begin{theorem}
\label{thm:1stable}
{\rm \cite{sa-1977}, \cite[Theorem~$9$]{eg-1984}}
A connected graph $G$ admits a $1$-stable trace if and only if $\delta(G) > 1$.
\end{theorem}

\begin{theorem}
\label{thm:2stable} 
{\rm \cite[Theorem~$3.1$]{kl-2013}}
A connected graph $G$ admits a $2$-stable trace if and only if $\delta(G) > 2$. 
\end{theorem}

Note that a vertex $v$ of degree $d$ implies that every double trace $W$ in $G$ has a repetition of order $d$. A leaf in a graph necessarily implies no double trace is $1$-stable, and similarly, a vertex $v$ of degree $2$ implies that a repetition of order $2$ is present in every double trace.  

Our key object in this paper is a \emph{strong trace}, which is a double trace without nontrivial repetitions. Observe that in a graph $G$ with $\delta(G) \ge 3$ every strong trace is $2$-stable. If also $\Delta(G) \le 5$ then every $2$-stable trace is also a strong one. 

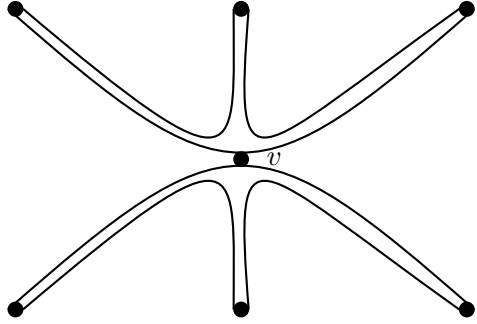
\begin{figure}[ht!]
\begin{center}
\begin{tikzpicture}[scale=1.0,style=thick]
\fill (3,2) circle (3pt);
\fill (3.2,2) node[right]{$v$};
\fill (0,0) circle (3pt);
\fill (3,0) circle (3pt);
\fill (6,0) circle (3pt);
\fill (0,4) circle (3pt);
\fill (3,4) circle (3pt);
\fill (6,4) circle (3pt);
\draw (-0.1,4)  .. controls (2.8,1.45) and (3.2,1.45) .. (6.1,4);
\draw (0.1,4)  .. controls (3.1,1.45) and (2.9,2) .. (2.9,4);
\draw (3.1,4)  .. controls (2.9,1.45) and (3.1,2) .. (5.9,4);
\draw (-0.1,0)  .. controls (2.8,2.55) and (3.2,2.55) .. (6.1,0);
\draw (0.1,0)  .. controls (3.1,2.55) and (2.9,2) .. (2.9,0);
\draw (3.1,0)  .. controls (2.9,2.55) and (3.1,2)  .. (5.9,0);
\end{tikzpicture}
\end{center}
\caption{A $3$-repetition in a vertex $v$ of degree $6$}
\label{fig:repetition}
\end{figure}

However, if $v$ is a vertex of degree at least $6$, then a stable trace $W$ may have a $3$-repetition at $v$, see Fig.~\ref{fig:repetition}. 

Our main results are the following theorems:
\begin{theorem}
\label{thm:strong}
Every connected graph $G$ admits a strong trace.
\end{theorem}

Now Theorem~\ref{thm:strong} implies:

\begin{theorem}
\label{thm:realize}
Every connected graph $G$ can be (at least in theory) constructed from a single coiled-coil-forming segment.
\end{theorem}

A weaker version of Theorem~\ref{thm:realize} (limited to graphs of polyhedra) was stated in~\cite{kl-2013}. The classical Steinitz' theorem~\cite{ste-1922} namely states that every 3-connected planar graph $G$ is isomorphic to a graph of a polyhedron --- $G=G(P)$ for some polyhedron $P$ --- and vice versa, for every polyhedron $P$ its graph $G(P)$ is a planar $3$-connected graph. Now $3$-connectivity of $G$ implies $\delta(G) \ge 3$, a condition heavily used in~\cite{kl-2013}.

We shall prove Theorem~\ref{thm:strong} in the next section after establishing the connection between strong traces and embeddings of graphs in surfaces.

%%%%%%%%%%%%%%%%%%%%%%%%%%%%%%%%%%%
\section{Graph embeddings and strong traces}
\label{sec:main}

In this section we establish the duality between embeddings of graphs in surfaces and strong traces in graphs. We shall first cover the necessary material on combinatorial embeddings of graphs. For more detail on the topic see~\cite{moh-2001}. 

A \emph{(combinatorial) embedding} of a graph $G$ in a surface $\Sigma$ is a pair $(\Pi, \lambda)$, where $\Pi = \{ \pi_v \mid v \in V(G)\}$ so that for every vertex $v \in V(G)$ $\pi_v$ is a \emph{cyclic permutation} of $E(v)$, and $\lambda : E(G) \rightarrow \{-1,1\}$. We shall call $\Pi$ the \emph{rotation system} and $\pi_v \in \Pi$ the \emph{local rotation around $v$}, whereas $\lambda$ is called the \emph{signature (of edges)}. 

The permutation $\pi_v$ describes the clock-wise ordering of edges emanating from $v$. For a pair of adjacent vertices $uv$ the signature $\lambda(uv)$ encodes the possible match of clockwise orientations $\pi_u$ and $\pi_v$: $\lambda(uv)=1$ if and only if the clockwise orientation around $u$ matches the one around $v$ when traversing the edge $uv$.   A \emph{facial walk} of $(\Pi, \lambda)$ is a closed walk in $G$ obtained by the following procedure. We start at an arbitrary vertex $u$, choose an arbitrary incident edge $uv$ and an initial signature value $\lambda_0 \in \{-1,1\}$. Now we repeat the following steps: move along the chosen edge $uv$, multiply the signature $\lambda_0$ by the signature of a traversed edge $\lambda(uv)$, and choose the next edge $vw$ so that either $\pi^{-1}_v(uv)=vw$ or $\pi_v(uv)=vw$ depending on whether $\lambda_0=1$ or $\lambda_0=-1$, respectively. We terminate the procedure when (i) we reach $u$, (ii) the next edge to travel is the initial edge $uv$, and (iii) $\lambda_0$ equals the initially chosen value.
We consider two facial walks the same if they only differ in respective initial vertices and/or their orientations.

The surface $\Sigma$ is uniquely determined by the combinatorial embedding $(\Pi, \lambda)$. $\Sigma$ is orientable if $G$ contains no cycle with an odd number of edges having negative signature, and is nonorientable if there exists a cycle $C$ having an odd number of edges with negative signature.  
The genus of $\Sigma$ is determined by the number of facial walks of $(\Pi, \lambda)$.

The sense of orientation changes at every edge with negative signature when traveling along $C$. If the number of changes along $C$ is odd a narrow strip around $C$ is homeomorphic to the M\"obius band. Face-wise --- by decreasing the number of facial walks we obtain surfaces of higher genera.

Two embeddings $(\Pi, \lambda)$ and $(\Pi', \lambda')$ are \emph{equivalent}  if one can be obtained from the other by repetitively replacing a single local rotation at $v$ by its inverse and at the same time altering signatures of every edge emanating from $v$.
 
Observe that every edge $uv$ of $G$ appears exactly twice in the collection of facial walks of $G$ in an embedding $(\Pi, \lambda)$.

An embedding $(\Pi, \lambda)$ determines the collection of facial walks, but it is also the other way around. A collection of closed walks ${\cal W}=\{W_1,\ldots,W_k\}$ so that every edge $uv \in E(G)$ 
appears exactly twice in $\cal W$ determines the embedding  $(\Pi, \lambda)$ up to equivalence: a sequence $e v e'$ along a facial walk implies that $e$ and $e'$ are consecutive in $\pi_v$, and a sequence $e v e' v' e''$ along a facial walk determines the signature of $e'$: $\lambda(e')=1$ if and only if either $\pi_v(e)=e'$ and $\pi_{v'}(e')=e''$ or $\pi_{v'}(e'')=e'$ and $\pi_v(e')=e$.

An alternative way to represent the surface $\Sigma$ is by taking its \emph{polygonal schema}: take a collection of disks, one per each facial walk in ${\cal W}=\{W_1,\ldots,W_k\}$, and make identification along their borders according to pairs of edges in $\cal W$. 

It is known that $\Sigma$ is orientable if the facial walks in $\cal W$ can be chosen in such a way that every edge is traversed twice in opposite directions. 

We proceed with a basic result on embeddings of connected graphs. A \emph{$k$-face embedding} is an embedding with exactly $k$ faces (facial walks). Next theorem was independently proven by Edmonds~\cite{ed-1965} and later Pisanski~\cite{pi-1978}.

\begin{theorem}
\label{thm:1face}
{\rm \cite{ed-1965}, \cite{pi-1978}}
 Every connected graph $G$ admits a $1$-face embedding in some surface $\Sigma$.
\end{theorem}

\proof
Let $(\Pi,\lambda)$ be a combinatorial embedding of $G$ with the smallest number of facial walks. If the number of facial walks is at least $2$, then some edge $e=uv$ is contained in a pair of distinct facial walks $W_1$ and $W_2$. We claim that changing the signature of $e$ reduces the number of facial walks by one. 

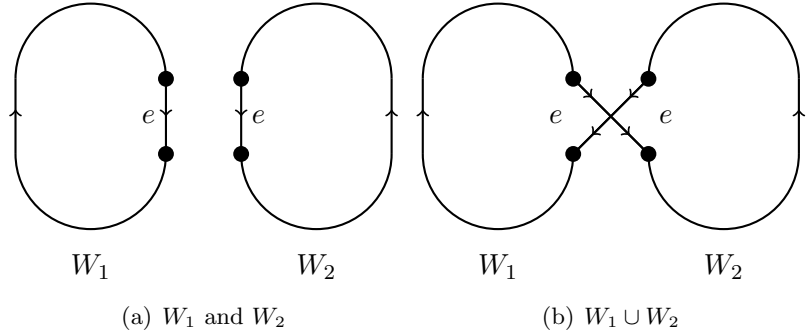
\begin{figure}[ht!]
\begin{center}
\subfigure[$W_1$ and $W_2$]
{
\begin{tikzpicture}[scale=1.0,style=thick]
\fill (1,3) circle (3pt);
\fill (1,2) circle (3pt);
\fill (2,3) circle (3pt);
\fill (2,2) circle (3pt);
\draw[decoration={markings, mark=at position 0.5 with {\arrow{>}}},postaction={decorate}] (1,3)--(1,2);
\draw[decoration={markings, mark=at position 0.5 with {\arrow{>}}},postaction={decorate}] (2,3)--(2,2);
\draw[decoration={markings, mark=at position 0.5 with {\arrow{<}}},postaction={decorate}] (-1,3)--(-1,2);
\draw[decoration={markings, mark=at position 0.5 with {\arrow{<}}},postaction={decorate}] (4,3)--(4,2);
\draw (1,3) arc (0:180:1);
\draw (-1,2) arc (-180:0:1);
\draw (4,3) arc (0:180:1);
\draw (2,2) arc (-180:0:1);
\draw (1,2.5) node[left] {$e$};
\draw (2,2.5) node[right] {$e$};
\draw (0,0.5) node {$W_1$};
\draw (3,0.5) node {$W_2$};
\end{tikzpicture}
}
\subfigure[$W_1 \cup W_2$]
{
\begin{tikzpicture}[scale=1.0,style=thick]
\fill (1,3) circle (3pt);
\fill (1,2) circle (3pt);
\fill (2,3) circle (3pt);
\fill (2,2) circle (3pt);
\draw[decoration={markings, mark=at position 0.25 with {\arrow{>}}},postaction={decorate}] (1,3)--(2,2);
\draw[decoration={markings, mark=at position 0.25 with {\arrow{>}}},postaction={decorate}] (2,3)--(1,2);
\draw[decoration={markings, mark=at position 0.75 with {\arrow{>}}},postaction={decorate}] (1,3)--(2,2);
\draw[decoration={markings, mark=at position 0.75 with {\arrow{>}}},postaction={decorate}] (2,3)--(1,2);
\draw[decoration={markings, mark=at position 0.5 with {\arrow{<}}},postaction={decorate}] (-1,3)--(-1,2);
\draw[decoration={markings, mark=at position 0.5 with {\arrow{<}}},postaction={decorate}] (4,3)--(4,2);
\draw (1,3) arc (0:180:1);
\draw (-1,2) arc (-180:0:1);
\draw (4,3) arc (0:180:1);
\draw (2,2) arc (-180:0:1);
\draw (1,2.5) node[left] {$e$};
\draw (2,2.5) node[right] {$e$};
\draw (0,0.5) node {$W_1$};
\draw (3,0.5) node {$W_2$};
\end{tikzpicture}
}
\end{center}
\caption{Construction from proof of Theorem~\ref{thm:1face}}
\label{fig:construction1}
\end{figure}

We may assume that $W_1$ and $W_2$ traverse $e$ in the same direction. Let us start walking along $W_1$. Continuing along $e$ with the change of its signature routes our walk following $W_2$, then back to $e$ where it continues along $W_1$. This implies that the walks $W_1$ and $W_2$ merge into a single facial walk in the adjusted embedding, see Fig.~\ref{fig:construction1}. The remaining facial walks clearly do not change.

\qed

An easy consequence of Theorem~\ref{thm:1face} is the classical theorem of Ringel.
\begin{theorem}
\label{thm:1face:non}
{\rm \cite[Theorem~$13$]{ri-1977}, \cite[Theorem~$8$]{st-1978}}
Every connected graph $G$ which is not a tree  has a $1$-face embedding in some nonorientable surface.
\end{theorem}

\proof
Assume that $G$ is not a tree and let $(\Pi,\lambda)$ be a $1$-face embedding of $G$ in some surface $\Sigma$. Such an embedding exists by Theorem~\ref{thm:1face}. 

Assume that $\Sigma$ is orientable, and let $W$ be the (only) facial walk which traverses every edge twice, once in every direction. As $G$ is not a tree there exists an edge $e=u_1v_1$ which is not a cutedge. We claim that changing the signature of $e$ produces a $1$-face embedding $(\Pi',\lambda')$ of $G$ into a nonorientable surface $\Sigma'$.

Let us denote 
$$W = u_0 \ldots f_1 u_1 e v_1 e_2 v_2 \ldots v_k e_k v_1 e u_1 g_1 \ldots u_0.$$
Altering the signature of $e$ yields an alternative embedding whose only facial walk equals
$$W' = u_0 \ldots f_1 u_1 e v_1 e_k v_k \ldots v_2 e_2 v_1 e u_1 g_1 \ldots u_0$$
obtained by reversing the subwalk between occurrences of $e$. 

As $e$ is not a cutedge there exists a cycle $C$ passing through $e$, which contains an odd number of edges whose $\lambda'$ signature is negative. Hence $\Sigma'$ is a nonorientable surface. 
\qed

Let 
$$W=v_0 e_1 v_1 \ldots v_{\ell-1} e_\ell v_\ell$$ 
be a double trace in $G$. Fix a vertex $v \in V(G)$ and let $E(v)$ be the set of edges emanating from $v$. Let us build a $2$-regular graph (a union of cycles) $F_{v,W}$, having $E(v)$ as its vertex set by making edges $e,e' \in E(v)$ adjacent if $e$ and $e'$ are consecutive edges along $W$ (where a $1$-repetition at $v$ constructs a loop and a $2$-repetition gives rise to a pair of parallel edges in $F_{v,W}$). The graph $F_{v,W}$ is also called the \emph{vertex figure of $v$} (with respect to a double trace $W$).

The connection between vertex figures and graph embeddings is best explained via the following proposition.

\begin{proposition}
\label{prp:strong:vfigure}
Let $G$ be a connected graph and $W$ a double trace in $G$. Then $W$ is a strong trace if and only if every vertex figure $F_{v,W}$ is a single cycle.
\end{proposition}

\proof
Assume first that $W$ is not strong. Then there exists a nontrivial $N$-repetition at some vertex $v \in V(G)$. Let $N' = N(v) \setminus N$, which is also nonempty. We claim that the vertex figure $F_{v,W}$ contains at least two cycles. 

Let $e$ be an edge incident with $v$ whose other endvertex lies in $N$. 
In the vertex figure $F_{v,W}$ the edge $e$ can only be adjacent to an edge from $E(v,N)$, and consequently none of the edges $e' \in E(v,N')$ lies in the same cycle of $F_{v,W}$ as $e$.

For the converse, let $W$ be a strong trace, and let us pick an arbitrary vertex $v$. Assume that $F_{v,W}$ contains a pair of disjoint cycles $C_1$ and $C_2$. Let $N$ be the set of endvertices of edges from $C_1$ different from $v$. Now $W$ contains an $N$-repetition at $v$, as entering $v$ from $N$ implies that $W$ also exits towards a vertex from $N$. As $N \ne \emptyset$ and $N \ne N(v)$ we have a nontrivial repetition at $v$ which is absurd.
\qed

Assume that $F_{v,W}$ is a single cycle. An orientation of $F_{v,W}$ can be interpreted as a cyclic permutation $\pi_v$ of $E(v)$. What is more, if every vertex figure is a single cycle, the collection of cyclic permutations $\Pi=\{\pi_v \mid v \in V(G)\}$ is the first component of an embedding, whose only facial cycle equals $W$. 

To sum it all up. By Theorem~\ref{thm:1face} every connected graph admits a $1$-face embedding $(\Pi,\lambda)$ into some closed surface $\Sigma$. The only facial cycle $W$ of this embedding is a double trace, and as every vertex figure $F_{v,W}$ is a single cycle, Proposition~\ref{prp:strong:vfigure} implies $W$ is strong. This completes the proof of Theorem~\ref{thm:strong}.

Theorem~\ref{thm:strong} easily implies the following proposition, which in turn implies both Theorem~\ref{thm:1stable} and Theorem~\ref{thm:2stable}.

\begin{proposition}
\label{prp:nstable}
Let $G$ be a connected graph. Then $G$ admits a $d$-stable trace if and only if $\delta(G)>d$.
\end{proposition}

\proof
It is enough to note that a strong trace in $G$ is $d$-stable, provided that no vertex in $G$ has degree $\le d$.
\qed

%%%%%%%%%%%%%%%%%%%%%%%%%%
\section{Antiparallel strong traces}
\label{sec:anti}

Let $W$ be a double trace in $G$. As mentioned in Section~\ref{sec:introd} every edge $e=uv$ of graph $G$ corresponds to a coiled-coil dimer and is thus traversed exactly twice in strong trace and $d$-stable trace $W$. If $W$ traverses $e$ in the same direction twice (either both times from $u$ to $v$ or both times from $v$ to $u$) then we call $e$ a {\em parallel edge} (with respect to $W$), otherwise $e$ is  an {\em antiparallel edge}. A double trace $W$ is a {\em parallel double trace} if every edge of $G$ is parallel and an {\em antiparallel double trace} if every edge of $G$ is antiparallel. 

The motivation for this concept also comes from self-assembly nanostructure design~\cite{gr-2013}. Parallel double traces represent polyhedra in which on every edge the two coiled-coil-forming segments would be aligned in the same direction while antiparallel double traces represent polyhedra in which on every edge two coiled-coil-forming segments would be aligned in the opposite direction. Because of apparent lack of polypeptide pairs which form antiparallel coiled-coil dimers~\cite{gr-2013}, especially detailed study of the first type would be of a great use. In this section we discuss antiparallel strong traces  and  turn to parallel strong traces in next section.

The main result of this section can be read as follows.

\begin{theorem}
\label{thm:santi}
A graph $G$ admits an antiparallel strong trace strong trace if and only if $G$ has a spanning tree $T$ such that each connected component of $G - E(T)$ has an even number of edges.
\end{theorem}

In the rest of this section we shall prove Theorem~\ref{thm:santi}.

Connection between antiparallel double traces and embeddings of graphs was (to some extent) already observed in~\cite{sk-1990} and~\cite{th-1990}. 

\begin{theorem}
\label{thm:anti-embedding}
A graph $G$ admits an antiparallel strong trace if and only if $G$ has an $1$-face embedding in some orientable surface.
\end{theorem}

\proof
Assume first that $G$ admits an antiparallel strong trace $W$. Now $W$ represents an unique facial walk of $G$ in an embedding $(\Pi, \lambda)$ and for every $v \in V(G)$ the vertex figure $F_{v,W}$ is a single cycle. Therefore $G$ has a $1$-face embedding in some surface $\Sigma$. Because every edge in $W$ is traversed twice in opposite direction, $\Sigma$ is orientable.
 
Conversely, let $(\Pi, \lambda)$ be a $1$-face embedding of $G$ in some orientable surface $\Sigma$. A $1$-face embedding $(\Pi, \lambda)$ determines an unique facial walk $W$. Clearly vertex figure $F_{v,W}$ is a single cycle for every $v \in V(G)$ and Proposition~\ref{prp:strong:vfigure} implies that $W$ is a strong trace in $G$. Because $\Sigma$ is orientable, every edge in $W$ is traversed twice in opposite directions, and is therefore antiparallel. 
\qed

Xuong characterized graphs which admit embeddings in orientable surface with at most 2 faces~\cite{xu-1979}. The {\em Betti number} of a graph $G$ is defined as $\beta(G) = |E(G)| - |V(G)| + 1$.
Observe also, as orientable surfaces have even Euler characteristics, that the number of faces in an orientable embedding of a graph $G$ is of different parity as its Betti number $\beta(G)$.

\begin{theorem}
\label{thm:upper}
{\rm \cite[Theorem~$2$]{xu-1979}}
A connected graph $G$ with even (odd) Betti number has an embedding in orientable surface with at most 2 faces if and only if it contains a spanning tree $T$ such that all (all but one) of connected components of $G - E(T)$ have an even number of edges.
\end{theorem}

A special case of Theorem~\ref{thm:upper} was later presented in~\cite{beh-1979} and~\cite{th-1995}: 

\begin{theorem}
\label{thm:strictly-upper}
A connected graph $G$ has an $1$-face embedding in an orientable surface if and only if $G$ has a spanning tree $T$ such that each connected component of $G - E(T)$ has an even number of edges.
\end{theorem}

To sum it all up. By Theorem~\ref{thm:anti-embedding}  a connected graph $G$ admits an antiparallel strong trace if and only if $G$ has an $1$-face embedding in some orientable surface. By Theorem~\ref{thm:strictly-upper} the latter is true if and only if $G$ has a spanning tree $T$ such that each connected component of $G - E(T)$ has an even number of edges. This completes the proof of Theorem~\ref{thm:santi}.

Already in $1895$ Tarry~\cite{tarry-1895} observed that every graph admits an antiparallel double trace. Almost a hundred years later Thomassen~\cite{th-1990} characterized graphs that admit antiparallel $1$-stable traces (thus solving a problem posed by Ore~\cite{ore-1951}):

\begin{theorem}
\label{thm:1anti}
{\rm \cite[Theorem~$3.3$]{th-1990}}
A graph $G$ admits an antiparallel $1$-stable trace if and only if $\delta(G) > 1$ and $G$ has a spanning tree $T$ such that each component of $G-E(T)$ either has an even number of edges or contains a vertex $v$ with $d_G(v)\ge 4$.
\end{theorem}

It would be interesting to characterize graphs which admit antiparallel $d$-stable traces, for every integer $d$.  By now it was observed that connection between graphs which admit antiparallel $d$-stable traces and pseudo-surfaces exists (for more on pseudo-surfaces, see~\cite{pot-2003}). Therefore the same approach as for characterization of graphs which admit antiparallel strong traces can not be used. Thus we pose:

\begin{problem}
Characterize graphs that admit an antiparallel $d$-stable trace for $d > 1$.
\end{problem}

%%%%%%%%%%%%%%%%%%%%%%%
\section{Parallel strong traces}
\label{sec:parallel}

We conclude with a characterization of graphs admitting parallel strong traces and parallel $d$-stable traces. Next proposition, which was observed in~\cite{kl-2013} easily follows if we traverse some Eulerian circuit of graph twice.

\begin{proposition}
\label{prp:parallel}
{\rm \cite[Proposition~$5.4$]{kl-2013}}
A graph $G$ admits a parallel $1$-stable trace if and only if $G$ is Eulerian. 
\end{proposition}

In~\cite{rus-2013} a similar theorem for parallel $2$-stable traces was proven.  Observe that in an Eulerian graph $G$ the condition $d_G(v) \ge 3$ is equivalent to $d_G(v) \ge 4$.

\begin{theorem}
\label{thm:parallel}
{\rm \cite[Theorem~$2.2$]{rus-2013}}
A graph $G$ admits a parallel $2$-stable trace if and only if $G$ is Eulerian and $\delta(G) \ge 3$. 
\end{theorem}

Our main result in this section is the following theorem, whose immediate corollary Theorem~\ref{thm:nparallel} easily implies both Proposition~\ref{prp:parallel} and Theorem~\ref{thm:parallel}.

\begin{theorem}
\label{thm:parallel:strong}
Let $G$ be a connected graph. $G$ admits a parallel strong trace if and only if $G$ is Eulerian.
\end{theorem}

\proof
If $G$ is not Eulerian then $G$ does not admit a parallel double trace --- the number of times a double trace enters a vertex $v$ of odd degree is on one hand odd (as it is equal to the number of times a double trace leaves) and even (as every edge incident with $v$ is either used twice or $0$ times for entering $v$), which is absurd.

For the converse direction let us consider a parallel double trace $W$ so that the collection of vertex figures $\{ F_{u,W} \mid u \in V(G)\}$ cumulatively has as few cycles as possible. If every vertex figure contains exactly one cycle, then by Proposition~\ref{prp:strong:vfigure} $W$ is a strong trace.

If on the other hand there exists a vertex whose vertex figure contains at least two cycles we shall be able to construct an alternative parallel double trace $W'$, so that the collection of alternative vertex figures $\{F_{u,W'} \mid u \in V(G)\}$ contains strictly fewer cycles. This will be the final contradiction in the proof.

Let $v$ be a vertex so that its vertex figure $F_{v,W}$ splits $E(v)$ into (at least) two cycles $C_1$ and  $C_2$. 
Choose an edge $e_1 = u_1v \in C_1$ so that $W$ uses $e_1$ in the direction towards $v$. Let $e_2=vu_2$ and $e_3=vu_3$ be the edges from $C_1$ that immediately succeed both occurrences  of $e_1$ along $W$ (note that $e_2$ may be equal to $e_3$). Next let $f_4=u_4v$ and $f_5=vu_5$ be edges from $C_2$ so that $u_4 f_4 v f_5 u_5$ is a subsequence of $W$.

\begin{figure}[ht!]
\begin{center}
\subfigure[$W$]
{
\begin{tikzpicture}[scale=1.0,style=thick]
\draw[dashed] (-2,5)--(0,5);
\draw[dashed] (-2,2)--(0,2);
\draw[dashed] (0,3)--(2,3);
\fill (-1,6) circle (3pt) node[left]{$u_1$};
\draw (-1,5.5) node[left] {$e_1$};
\fill (-1,5) circle (3pt) node[above right]{$v$};
\draw (-1,4.5) node[left] {$e_2$};
\fill (-1,4) circle (3pt) node[left]{$u_2$};
\draw (-1,3.5) node[left] {$A$};
\fill (-1,3) circle (3pt) node[left]{$u_1$};
\draw (-1,2.5) node[left] {$e_1$};
\fill (-1,2) circle (3pt) node[above right]{$v$};
\draw (-1,1.5) node[left] {$e_3$};
\fill (-1,1) circle (3pt) node[left]{$u_3$};

\draw (1,1) node[right] {$B$};
\fill (1,4) circle (3pt) node[right]{$u_5$};
\draw (1,3.5) node[right] {$f_5$};
\fill (1,3) circle (3pt) node[above left]{$v$};
\draw (1,2.5) node[right] {$f_4$};
\fill (1,2) circle (3pt) node[right]{$u_4$};
\draw (1,6) node[right] {$C$};

\draw[decoration={markings, mark=at position 0.5 with {\arrow{>}}},postaction={decorate}] (-1,6)--(-1,5);
\draw[decoration={markings, mark=at position 0.5 with {\arrow{>}}},postaction={decorate}] (-1,5)--(-1,4);
\draw[decoration={markings, mark=at position 0.5 with {\arrow{>}}},postaction={decorate}] (-1,4)--(-1,3);
\draw[decoration={markings, mark=at position 0.5 with {\arrow{>}}},postaction={decorate}] (-1,3)--(-1,2);
\draw[decoration={markings, mark=at position 0.5 with {\arrow{>}}},postaction={decorate}] (-1,2)--(-1,1);
\draw (1,1)--(1,2);
\draw[decoration={markings, mark=at position 0.5 with {\arrow{>}}},postaction={decorate}] (1,2)--(1,3);
\draw[decoration={markings, mark=at position 0.5 with {\arrow{>}}},postaction={decorate}] (1,3)--(1,4);
\draw (1,4)--(1,6);
\draw[decoration={markings, mark=at position 0.75 with {\arrow{>}}},postaction={decorate}] (-1,1) arc (-180:0:1);
\draw[decoration={markings, mark=at position 0 with {\arrow{>}}},postaction={decorate}] (1,6) arc (0:180:1);
\end{tikzpicture}
}
\subfigure[$W'$]
{
\begin{tikzpicture}[scale=1.0,style=thick]
\draw[dashed] (-2,5)--(0,5);
\draw[dashed] (-2,2)--(0,2);
\draw[dashed] (0,3)--(2,3);
\fill (-1,6) circle (3pt) node[left]{$u_1$};
\draw (-1,5.5) node[left] {$e_1$};
\fill (-1,5) circle (3pt) node[above right]{$v$};
\draw (-1,4.5) node[left] {$e_3$};
\fill (-1,4) circle (3pt) node[left]{$u_3$};
\draw (-1,3.5) node[left] {$B$};
\fill (-1,3) circle (3pt) node[left]{$u_4$};
\draw (-1,2.5) node[left] {$f_4$};
\fill (-1,2) circle (3pt) node[above right]{$v$};
\draw (-1,1.5) node[left] {$e_2$};
\fill (-1,1) circle (3pt) node[left]{$u_2$};

\draw (1,1) node[right] {$A$};
\fill (1,4) circle (3pt) node[right]{$u_5$};
\draw (1,3.5) node[right] {$f_5$};
\fill (1,3) circle (3pt) node[above left]{$v$};
\draw (1,2.5) node[right] {$e_1$};
\fill (1,2) circle (3pt) node[right]{$u_1$};
\draw (1,6) node[right] {$C$};

\draw[decoration={markings, mark=at position 0.5 with {\arrow{>}}},postaction={decorate}] (-1,6)--(-1,5);
\draw[decoration={markings, mark=at position 0.5 with {\arrow{>}}},postaction={decorate}] (-1,5)--(-1,4);
\draw[decoration={markings, mark=at position 0.5 with {\arrow{>}}},postaction={decorate}] (-1,4)--(-1,3);
\draw[decoration={markings, mark=at position 0.5 with {\arrow{>}}},postaction={decorate}] (-1,3)--(-1,2);
\draw[decoration={markings, mark=at position 0.5 with {\arrow{>}}},postaction={decorate}] (-1,2)--(-1,1);
\draw (1,1)--(1,2);
\draw[decoration={markings, mark=at position 0.5 with {\arrow{>}}},postaction={decorate}] (1,2)--(1,3);
\draw[decoration={markings, mark=at position 0.5 with {\arrow{>}}},postaction={decorate}] (1,3)--(1,4);
\draw (1,4)--(1,6);
\draw[decoration={markings, mark=at position 0.75 with {\arrow{>}}},postaction={decorate}] (-1,1) arc (-180:0:1);
\draw[decoration={markings, mark=at position 0 with {\arrow{>}}},postaction={decorate}] (1,6) arc (0:180:1);
\end{tikzpicture}
}
\end{center}
\caption{Construction from proof of Theorem~\ref{thm:parallel:strong}}
\label{fig:construction2}
\end{figure}
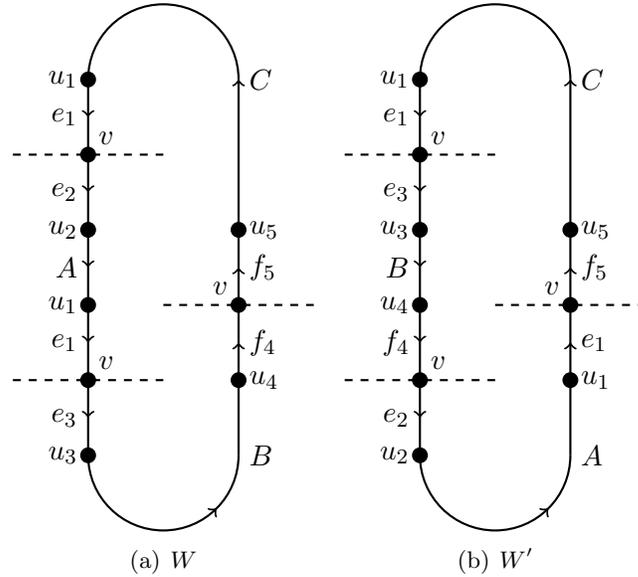

Without loss of generality (by choosing an alternative initial vertex along $W$) we may assume that
$$W = \ldots u_1 e_1 v e_2 u_2 \ldots u_1 e_1 v e_3 u_3 \ldots u_4 f_4 v f_5 u_5 \ldots .$$
Observe the following walk
$$W'= \ldots u_1 e_1 v e_3 u_3 \ldots u_4 f_4 v e_2 u_2 \ldots u_1 e_1 v f_5 u_5 \ldots$$
obtained by interchanging the two ``interior'' subwalks of $W$ between the three shown occurrences of $v$ in $W$, also see Fig.~\ref{fig:construction2}.

As $W'$ traverses the same collection of edges (in the same direction) as $W$, the walk $W'$ is indeed a parallel double trace. If $u \ne v$ then the new vertex figure $F_{u,W'}$ equals the original vertex figure $F_{u,W}$, since every pair $e,e'$ of edges meeting at $u$ are consecutive along $W'$ if and only if they are consecutive along $W$.

Now $W'$ only changes pairs of consecutive edges from $C_1 \cup C_2$, hence the only possible cycles of $F_{v,W'}$ which are not present in $F_{v,W}$ consist of edges from $C_1 \cup C_2$. Now the adjacencies $e_1-e_2$ and $f_4-f_5$ were replaced by $e_2-f_4$ and $e_1-f_5$ which implies that $C_1$ and $C_2$ merge into exactly one new cycle in $F_{v,W'}$ containing all edges from $C_1 \cup C_2$. Hence the total number of cycles in vertex figures has decreased by exactly one, which concludes the proof. 
\qed

The next theorem easily follows:

\begin{theorem}
\label{thm:nparallel}
A connected graph $G$ admits a parallel $d$-stable trace if and only if $G$ is Eulerian and $\delta(G) > d$.
\end{theorem}

Note that the construction used in~\cite{rus-2013} for the proof of  Theorem~\ref{thm:parallel} could be extended to yield an alternate proof of Theorem~\ref{thm:parallel:strong}.

%%%%%%%%%%%%%%%%%
\section{Conclusion}

Let us finish with a pair of open problems. We have provided a model for constructing a polypeptide nanostructure using a strong trace in the corresponding graph. A strong trace is a particular version of a closed walk, but is nevertheless encoded by a sequence which pins out its initial and terminal vertex. 
\begin{itemize}
\item Should we care which vertex of a nanostructure should be the initial vertex of an encoding of a double trace? How do both physical and chemical properties of the structure relate to the initial node of the polypeptide chain. 
\end{itemize}  
A fixed graph $G$ contains many strong traces, \cite{gr-2013} quotes 40 strong traces for the cube graph, for example.
Putting the initial vertex and the orientation aside, there still are many options and criteria along which one may choose a supposedly better strong trace.  
\begin{itemize}
\item Given a nanostructure, which strong trace to choose in order to maximize the probability an appropriate polypeptide chain will self assemble into the desired structure. 
\end{itemize}

%%%%%%%%%%%%%%%%%%%%%%%
\section*{Acknowledgements}
We are grateful to Michal Kotrb\v cik and Thomas W. Tucker for several remarks and suggestions which were of great help. Work supported in part by the ARSS of Slovenia, Research Grants
P1-0294, P1-0297  and N1-0011: GReGAS, supported in part by the European Science Foundation.

\baselineskip=0.21in

%%%%%%%%%%%%%%%%%%%%%

\end{document}